\documentclass[letterpaper, 10pt, journal]{IEEEtran}
\usepackage{cite}
\usepackage{amsmath}
\usepackage{amsfonts}
\usepackage{amsthm}
\theoremstyle{definition}

\usepackage{bbm}
\usepackage{empheq, bm, physics}
\usepackage{algorithm}
\usepackage{algorithmicx}
\usepackage{algpseudocode}
\usepackage{tikz}
\usepackage{dsfont}

\usepackage{tikzit}

\tikzstyle{new style 0}=[fill=white, draw=black, shape=circle]
\tikzstyle{new style 1}=[fill=white, draw=black, shape=rectangle]

\tikzstyle{new edge style 0}=[->, ultra thick, dashed]
\tikzstyle{new edge style 1}=[-, ultra thick, dashed, draw={rgb,255: red,21; green,115; blue,255}]
\tikzstyle{right}=[->, ultra thick]
\tikzstyle{left}=[<-, ultra thick]

\usepackage{float}

\usepackage{array}
\usepackage{physics}
\newtheorem{remark}{Remark}

\renewcommand{\skew}[1]{ \mathrm{ Skew } ( #1 ) }
\newcommand{\sym}[1]{\mathrm{Sym}_+( #1 )}
\newcommand{\sk}[1]{\mathrm{Skew}( #1 )}

\newcommand{\lag}{\mathcal{L}_\rho}
\newcommand{\blag}{\bar{\mathcal{L}}_\rho}
\newcommand{\lagr}[1]{\mathcal{L}_{\rho_{ #1 }}\qty}

\newcommand{\rgrad}{\mathrm{grad}\, }
\newcommand{\dist}[1]{\mathrm{dist}^{( #1 )}\, }

\makeatletter
\newcommand\fs@spaceruled{\def\@fs@cfont{\bfseries}\let\@fs@capt\floatc@ruled
  \def\@fs@pre{\vspace{\baselineskip}\hrule height.8pt depth0pt \kern2pt}%
  \def\@fs@post{\kern2pt\hrule\relax}%
  \def\@fs@mid{\kern2pt\hrule\kern2pt}%
  \let\@fs@iftopcapt\iftrue}
\makeatother

\begin{document}
\bstctlcite{IEEEexample:BSTcontrol}
%
\title{$H^2$-Optimal Reduction of Positive Networks using Riemannian Augmented Lagrangian Method}
%
%
%
\newcommand{\idnz}{\mathrm{nn}\, }
\newcommand{\idz}{\mathrm{z}\, }

\author{Sota~Misawa
        and~Kazuhiro~Sato
\thanks{S. Misawa and K. Sato are
with the Graduate School of Information Science and Technology, The University of Tokyo, Tokyo 113-8656, Japan, email:misawa-algebra@g.ecc.u-tokyo.ac.jp (S. Misawa), kazuhiro@mist.i.u-tokyo.ac.jp (K. Sato)}
}

\maketitle
\thispagestyle{empty} 
\pagestyle{empty}    
\begin{abstract}
In this study, we formulate the model reduction problem of a stable and positive network system as a constrained Riemannian optimization problem with the $H^2$-error objective function of the original and reduced network systems.
We improve the reduction performance of the clustering-based method, which is one of the most known methods for model reduction of positive network systems, by using the output of the clustering-based method as the initial point for the proposed method.
The proposed method reduces the dimension of the network system while preserving the properties of stability, positivity, and interconnection structure by applying the Riemannian augmented Lagrangian method (RALM) and deriving the Riemannian gradient of the Lagrangian.
To check the efficiency of our method, we conduct a numerical experiment and compare it with the clustering-based method in the sense of $H^2$-error and $H^\infty$-error.
\end{abstract}


%
\IEEEpeerreviewmaketitle

\section{Introduction}
%
%
%
%
Model reduction is a crucial step in designing controllers for large-scale network systems,
and thus some reduction methods have been proposed such as the balanced truncation (BT) method \cite{moore1981principal} and the application of the iterative rational Krylov algorithm (IRKA)  \cite{gugercin2008h_2, kellems2009low, gugercin2012structure, baur2011interpolatory, benner2015survey}. 
However, the BT and IRKA do not preserve the original interconnection structure in spite of the importance for
controlling and monitoring \cite{summers2014optimal, gates2016control, kim2018role, ishizaki2018graph}.
To resolve this issue,
the interconnection structure preserving model reduction methods for network systems have been proposed for a few decades.
For example,
 \cite{martin2018large} proposed to preserve the scale-free property of networks by formulating the interconnection constraints as the eigenvector centrality. In this method, the reduced network remains a flow network if the initial network is a flow network. Furthermore, \cite{vandendorpe2008model} introduced network reduction methods which preserve the interconnection structure of subsystems.

Moreover, model reduction methods of positive network systems whose outputs are always nonnegative under nonnegative inputs are important, because
the systems are often found in real world applications such as
pharmacokinetics, metabolism, epidemiology, ecology, and logistics \cite{haddad2010nonnegative, jacquez1993qualitative}.
In a positivity-preserving manner, some model reduction methods have been proposed in \cite{Cheng2020ModelRM, ISHIZAKI2015238, reis2009positivity}.
For instance, the clustering-based method \cite{ISHIZAKI2015238, reis2009positivity} is one of the most known methods which preserves the positivity.
Furthermore, in \cite{ISHIZAKI2015238, reis2009positivity}, the theoretical bounds of the $H^2$-error between the original and reduced systems are provided.
Nonetheless, the clustering-based method does not guarantee the $H^2$-optimality between the original and reduced network systems.

The $H^2$ optimal model reduction methods have been proposed in \cite{sato2015riemannian, sato2017structure, sato2018riemannian, sato2019riemannian}
based on Riemannian optimization \cite{absil2009optimization}.
In particular, in \cite{sato2019riemannian}, the set of stable matrices was endowed with the geometry of a Riemannian manifold, as summarized in Section II. D in this paper.
That is, the Riemannian optimization algorithm in \cite{sato2019riemannian} always produces a reduced asymptotically stable system at each iteration.
However, the above methods are not appropriate for network systems, because
the resulting reduced systems do not preserve the original interconnection structure.
That is, it is difficult to physically interpret the resulting reduced model.

For the $H^2$ optimal network system reduction, \cite{cheng2019reduced} proposed the $H^2$ optimal reduction method for linear consensus networks consisting of diffusively coupled single-integrators, which uses the clustering-based model reduction method as an initial network of the algorithm and aims to minimize the $H^2$ error by selecting suitable edge weights of the reduced network. Moreover, the method preserves the interconnection structure of the original network system. 
However, since this method is based on matrix inequalities, it is difficult to simply extend to general positive networks to preserve the positive property.

Therefore, to reduce large-scale asymptotically stable positive network systems, we formulate a novel $H^2$ optimization problem as a Riemannian optimization problem with constraints.
The introduction of the constraints is for preserving the positivity and interconnection structure of the original network system, and is the major difference with the existing problem formulations in \cite{sato2015riemannian, sato2017structure, sato2018riemannian, sato2019riemannian}.
To define the constraints, we use the result of a clustering-based model reduction method.
That is, the problem formulation in this paper can be regarded as a generalization of that in \cite{sato2019riemannian}.
 The main contribution is to develop the Riemannian augmented Lagrangian method (RALM) \cite{liu2019simple} for solving the problem.
 That is, the RALM
    preserves not only the stability and positivity but also the interconnection structure of the original system.
 To this end, we derive the Riemannian gradients, that are different from those of the objective function in \cite{sato2019riemannian}, of the Lagrangian composed of the objective function and a penalty term.

 The remainder of this paper is organized as follows.
In Section II, we describe the preliminary knowledge about asymptotically stable positive network systems, clustering-based network reduction, and the Riemannian manifold of stable matrices.
 In Section III, we formulate a novel $H^2$ optimal model reduction problem for preserving the positivity and interconnection structure of the original network system as a Riemannian optimization problem with constraints. 
 In Section IV, we propose the RALM algorithm by deriving the gradients of the Lagrangian.
In Section V, to illustrate the effectiveness of our method, we conduct a numerical experiment on a stable positive network and compare it with the clustering-based method from the viewpoints of $H^2$-error and $H^\infty$-error.
Finally, our conclusions are presented in Section VI.

\section{Preliminaries}
\subsection{Notation}
For a Riemannian manifold $M$, the tangent space at $x\in M$ is denoted by $T_xM$. We remark $\langle\cdot, \cdot\rangle_x^{(M)} : T_xM\times T_xM\rightarrow \mathbb{R}$ is an inner product at $x\in M$. 
For a vector $v\in \mathbb {R}^m$, $\|v\|_2$ denotes the usual Euclidean norm.
The $L^2$ space on $\mathbb{R}^m$ for $m\in \mathbb{Z}_{>0}$ is denoted by $L^2({\mathbb R}^m)$ with the norm $\|f\|_{L^2}:= \sqrt{\int_0^{\infty} \|f(t)\|_2^2 dt}$, where $f:[0,\infty)\rightarrow \mathbb {R}^m$ is a measurable function.
The $H^2$ and $H^{\infty}$ norms of a linear system whose transfer function is $G$ are defined by $\|G\|_{H^2} = \sqrt{\frac{1}{2\pi}\int_{-\infty}^\infty \norm{G(i\omega)}^2_\mathrm{F}\dd \omega}$ and $\|G\|_{H^{\infty}}=\sup_{\omega\in\mathbb{R}}\sigma_{\mathrm{max}}(G(i\omega))$, respectively,
where $\|\cdot\|_{\mathrm{F}}$ is the Frobenius norm and $\sigma_{\mathrm{max}}(G(i\omega))$ denotes the maximum singular value of $G(i\omega)$.

\subsection{Asymptotically Stable Positive Network System}
In this paper, we consider
\begin{align}
\begin{cases}
\dot{x}(t) = Ax(t) + Bu(t), \\
y(t) = Cx(t)
\end{cases} \label{original}
\end{align}
as the original large-scale network system
with the state $x(t)\in \mathbb{R}^n$, input $u(t)\in \mathbb{R}^m$, output $y(t)\in \mathbb{R}^p$, and coefficient matrices $A\in \mathbb{R}^{n\times n}$, $B\in \mathbb{R}^{n\times m}$, and $C\in \mathbb{R}^{p\times n}$.

We assume the following:
\begin{enumerate}
    \item The original network system (\ref{original}) is asymptotically stable.
That is,
 the real parts of all the eigenvalues of the matrix $A$ are negative. 
In this case, we call $A$ a stable matrix.

The matrix $A$ is Metzler.
That is, 
  every off-diagonal entry of $A$ is nonnegative.
  
\item The matrices $B$ and $C$ are nonnegative.
\end{enumerate}

Assumptions 2) and 3) mean that the original network system (\ref{original}) is essentially positive, as shown in \cite{haddad2010nonnegative}.
That is, not only the output $y(t)$ but also the state $x(t)$ is nonnegative under the nonnegative input $u(t)$ and initial state $x(0)$.
In fact, the solution to system \eqref{original} is given by
    $x(t) = \exp (At) x(0) + \int_0^t \exp (A(t-\tau))Bu(\tau) {\rm d}\tau$. 
If $A$ is Metzler, $\exp(At)$ for any $t\in \mathbb{R}$ is nonnegative.
Thus, 
Assumptions 2) and 3) imply that $x(t)$ and $y(t)$ are nonnegative under the nonnegative input $u(t)$ and initial state $x(0)$.

We denote the original network graph by $\mathcal{G}=(\mathcal{V}, \mathcal{E})$.

\subsection{Network Reduction Based on Clustering}

For fixed $r$ ($<n$), we reduce the original system (\ref{original}) to a $r$-dimensional system
\begin{align}
  \label{eq:reduced system all}
  \begin{cases}
    \dot{x}_r (t) = A_r x_r (t) + B_r u(t)\\
    y_r (t) = C_r x_r (t), 
  \end{cases}
\end{align}
where $x_r(t)\in \mathbb{R}^r$, $y_r(t)\in \mathbb{R}^p$, 
and $(A_r, B_r, C_r) \in \mathbb{R}^{r\times r}\times \mathbb{R}^{r\times m}\times \mathbb{R}^{p\times r}$. 
We define
$\mathcal{G}_r=(\mathcal{V}_r, \mathcal{E}_r)$
as the reduced network graph associated to system \eqref{eq:reduced system all}.

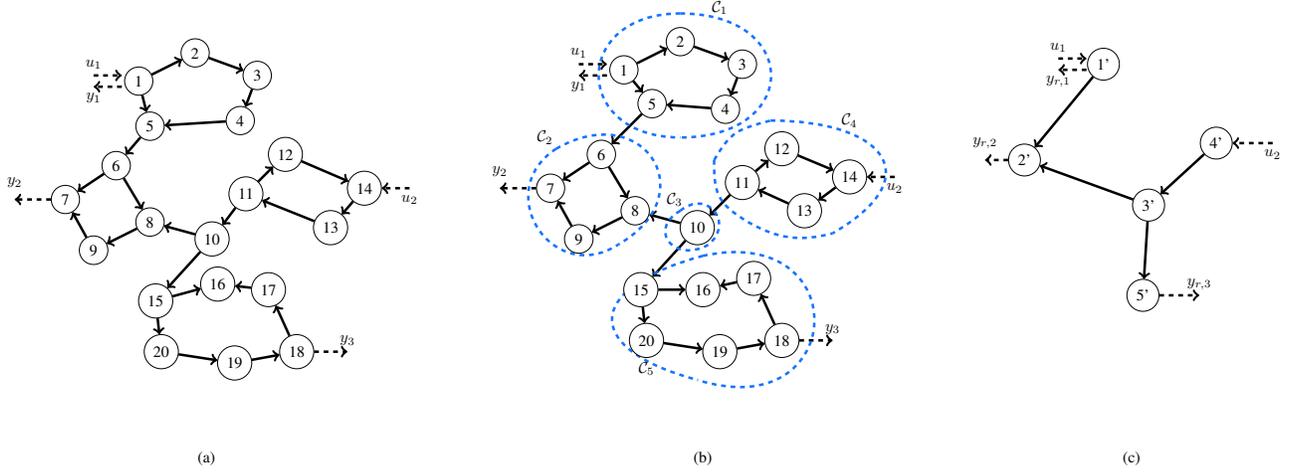
\begin{figure*}[t]
\vspace{3mm}
\begin{center}
\scalebox{0.6}{\begin{tikzpicture}
	\begin{pgfonlayer}{nodelayer}
		\node [style=new style 0] (1) at (-25, 4.75) {1};
		\node [style=new style 0] (2) at (-19.75, 5) {3};
		\node [style=new style 0] (3) at (-22.5, 6) {2};
		\node [style=new style 0] (4) at (-20.5, 3) {4};
		\node [style=new style 0] (5) at (-24.5, 2.75) {5};
		\node [style=new style 0] (6) at (-26, 1) {6};
		\node [style=new style 0] (7) at (-28.25, -0.5) {7};
		\node [style=new style 0] (8) at (-27, -2.75) {9};
		\node [style=new style 0] (9) at (-24.5, -1.5) {8};
		\node [style=new style 0] (10) at (-21.75, -2.25) {10};
		\node [style=new style 0] (11) at (-18.5, 1.5) {12};
		\node [style=new style 0] (12) at (-20.25, -0.25) {11};
		\node [style=new style 0] (13) at (-16.5, -1.75) {13};
		\node [style=new style 0] (14) at (-15, 0) {14};
		\node [style=new style 0] (15) at (-24.25, -5) {15};
		\node [style=new style 0] (16) at (-21.5, -4.25) {16};
		\node [style=new style 0] (17) at (-19.25, -4.5) {17};
		\node [style=new style 0] (18) at (-24, -7.25) {20};
		\node [style=new style 0] (19) at (-20.75, -7.75) {19};
		\node [style=new style 0] (20) at (-18, -7.25) {18};
		\node [style=none] (21) at (-27, 4) {$y_1$};
		\node [style=none] (22) at (-13, -0.5) {$u_2$};
		\node [style=none] (23) at (-15.75, -6.75) {$y_3$};
		\node [style=none] (24) at (-30.5, 0.25) {$y_2$};
		\node [style=none] (25) at (-27, 5.5) {$u_1$};
		\node [style=none] (26) at (-27, 5) {};
		\node [style=none] (27) at (-13, 0) {};
		\node [style=none] (28) at (-15.75, -7.25) {};
		\node [style=none] (29) at (-30.5, -0.5) {};
		\node [style=none] (30) at (-27, 4.5) {};
		\node [style=none] (31) at (-25.75, 5) {};
		\node [style=none] (32) at (-25.75, 4.5) {};
		\node [style=new style 0] (33) at (-3.5, 5.25) {1};
		\node [style=new style 0] (34) at (1.75, 5.5) {3};
		\node [style=new style 0] (35) at (-1, 6.5) {2};
		\node [style=new style 0] (36) at (1, 3.5) {4};
		\node [style=new style 0] (37) at (-2.25, 3.75) {5};
		\node [style=new style 0] (38) at (-4.5, 1.5) {6};
		\node [style=new style 0] (39) at (-6.75, 0) {7};
		\node [style=new style 0] (40) at (-5.5, -2.25) {9};
		\node [style=new style 0] (41) at (-3, -1) {8};
		\node [style=new style 0] (42) at (-0.25, -1.75) {10};
		\node [style=new style 0] (43) at (3.5, 1.75) {12};
		\node [style=new style 0] (44) at (1.75, 0.25) {11};
		\node [style=new style 0] (45) at (4.5, -1) {13};
		\node [style=new style 0] (46) at (6.5, 0.5) {14};
		\node [style=new style 0] (47) at (-2.75, -4.5) {15};
		\node [style=new style 0] (48) at (0, -4.5) {16};
		\node [style=new style 0] (49) at (2.25, -4) {17};
		\node [style=new style 0] (50) at (-2.5, -6.75) {20};
		\node [style=new style 0] (51) at (0.75, -7.25) {19};
		\node [style=new style 0] (52) at (3.5, -6.75) {18};
		\node [style=none] (53) at (-5.5, 4.5) {$y_1$};
		\node [style=none] (54) at (8.5, 0) {$u_2$};
		\node [style=none] (55) at (5.75, -6.25) {$y_3$};
		\node [style=none] (56) at (-9, 0.75) {$y_2$};
		\node [style=none] (57) at (-5.5, 6) {$u_1$};
		\node [style=none] (58) at (-5.5, 5.5) {};
		\node [style=none] (59) at (8.5, 0.5) {};
		\node [style=none] (60) at (5.75, -6.75) {};
		\node [style=none] (61) at (-9, 0) {};
		\node [style=none] (62) at (-5.5, 5) {};
		\node [style=none] (63) at (-4.25, 5.5) {};
		\node [style=none] (64) at (-4.25, 5) {};
		\node [style=none] (65) at (3, 3) {};
		\node [style=none] (66) at (5.75, -2) {};
		\node [style=none] (67) at (-1, 7.75) {};
		\node [style=none] (68) at (-1, 2.25) {};
		\node [style=none] (69) at (-6.5, 1.75) {};
		\node [style=none] (70) at (-3.75, -2.5) {};
		\node [style=none] (71) at (0, -3) {};
		\node [style=none] (72) at (-0.25, -8.5) {};
		\node [style=none] (73) at (0.75, 8) {};
		\node [style=none] (74) at (0.75, 8) {$\mathcal{C}_1$};
		\node [style=none] (76) at (-7, 2.25) {$\mathcal{C}_2$};
		\node [style=none] (77) at (4.25, -8.75) {};
		\node [style=none] (78) at (-2.5, -8) {$\mathcal{C}_5$};
		\node [style=none] (79) at (6.5, 3) {$\mathcal{C}_4$};
		\node [style=none] (80) at (-1, -1) {};
		\node [style=none] (81) at (0.25, -2.5) {};
		\node [style=none] (82) at (-1.25, -0.5) {$\mathcal{C}_3$};
		\node [style=none] (83) at (2, -1.5) {};
		\node [style=none] (84) at (0.5, 1.25) {};
		\node [style=new style 0] (85) at (17.75, 5.5) {1'};
		\node [style=new style 0] (86) at (14.25, 1.25) {2'};
		\node [style=new style 0] (87) at (19.75, -0.75) {3'};
		\node [style=new style 0] (88) at (22.75, 2) {4'};
		\node [style=new style 0] (89) at (19.5, -4.75) {5'};
		\node [style=none] (90) at (15.75, 4.75) {$y_{r, 1}$};
		\node [style=none] (91) at (15.75, 6.25) {$u_1$};
		\node [style=none] (92) at (15.75, 5.75) {};
		\node [style=none] (93) at (15.75, 5.25) {};
		\node [style=none] (94) at (17, 5.75) {};
		\node [style=none] (95) at (17, 5.25) {};
		\node [style=none] (96) at (12.5, 2) {$y_{r, 2}$};
		\node [style=none] (97) at (12.5, 1.25) {};
		\node [style=none] (98) at (25.25, 1.5) {$u_2$};
		\node [style=none] (99) at (25.25, 2) {};
		\node [style=none] (100) at (22, -4.25) {$y_{r, 3}$};
		\node [style=none] (101) at (22, -4.75) {};
		\node [style=none] (102) at (-22, -12) {(a)};
		\node [style=none] (103) at (0, -12) {(b)};
		\node [style=none] (104) at (19, -12) {(c)};
	\end{pgfonlayer}
	\begin{pgfonlayer}{edgelayer}
		\draw [style=right] (1) to (3);
		\draw [style=right] (3) to (2);
		\draw [style=right] (2) to (4);
		\draw [style=right] (4) to (5);
		\draw [style=left] (5) to (1);
		\draw [style=right] (5) to (6);
		\draw [style=right] (6) to (7);
		\draw [style=right] (8) to (7);
		\draw [style=left] (8) to (9);
		\draw [style=left] (9) to (6);
		\draw [style=left] (9) to (10);
		\draw [style=right] (12) to (10);
		\draw [style=right] (12) to (11);
		\draw [style=right] (11) to (14);
		\draw [style=right] (14) to (13);
		\draw [style=right] (13) to (12);
		\draw [style=right] (10) to (15);
		\draw [style=right] (15) to (16);
		\draw [style=left] (16) to (17);
		\draw [style=left] (17) to (20);
		\draw [style=left] (20) to (19);
		\draw [style=left] (19) to (18);
		\draw [style=left] (18) to (15);
		\draw [style=new edge style 0] (26.center) to (31.center);
		\draw [style=new edge style 0] (32.center) to (30.center);
		\draw [style=new edge style 0] (27.center) to (14);
		\draw [style=new edge style 0] (20) to (28.center);
		\draw [style=new edge style 0] (7) to (29.center);
		\draw [style=right] (33) to (35);
		\draw [style=right] (35) to (34);
		\draw [style=right] (34) to (36);
		\draw [style=right] (36) to (37);
		\draw [style=left] (37) to (33);
		\draw [style=right] (37) to (38);
		\draw [style=right] (38) to (39);
		\draw [style=right] (40) to (39);
		\draw [style=left] (40) to (41);
		\draw [style=left] (41) to (38);
		\draw [style=left] (41) to (42);
		\draw [style=right] (44) to (42);
		\draw [style=right] (44) to (43);
		\draw [style=right] (43) to (46);
		\draw [style=right] (46) to (45);
		\draw [style=right] (45) to (44);
		\draw [style=right] (42) to (47);
		\draw [style=right] (47) to (48);
		\draw [style=left] (48) to (49);
		\draw [style=left] (49) to (52);
		\draw [style=left] (52) to (51);
		\draw [style=left] (51) to (50);
		\draw [style=left] (50) to (47);
		\draw [style=new edge style 0] (58.center) to (63.center);
		\draw [style=new edge style 0] (64.center) to (62.center);
		\draw [style=new edge style 0] (59.center) to (46);
		\draw [style=new edge style 0] (52) to (60.center);
		\draw [style=new edge style 0] (39) to (61.center);
		\draw [style=new edge style 1, in=30, out=0, looseness=2.25] (65.center) to (66.center);
		\draw [style=new edge style 1, in=0, out=0, looseness=2.50] (67.center) to (68.center);
		\draw [style=new edge style 1, in=180, out=180, looseness=2.25] (67.center) to (68.center);
		\draw [style=new edge style 1, bend right=90, looseness=1.75] (69.center) to (70.center);
		\draw [style=new edge style 1, in=30, out=30, looseness=2.25] (69.center) to (70.center);
		\draw [style=new edge style 1, bend left=105, looseness=3.25] (71.center) to (72.center);
		\draw [style=new edge style 1, in=165, out=-165, looseness=2.50] (71.center) to (72.center);
		\draw [style=new edge style 1, bend left=90, looseness=2.00] (80.center) to (81.center);
		\draw [style=new edge style 1, in=210, out=-150, looseness=2.25] (80.center) to (81.center);
		\draw [style=new edge style 1, in=-165, out=60, looseness=1.25] (84.center) to (65.center);
		\draw [style=new edge style 1, bend right, looseness=1.25] (84.center) to (83.center);
		\draw [style=new edge style 1, in=195, out=-30] (83.center) to (66.center);
		\draw [style=right] (85) to (86);
		\draw [style=left] (86) to (87);
		\draw [style=left] (87) to (88);
		\draw [style=right] (87) to (89);
		\draw [style=new edge style 0] (92.center) to (94.center);
		\draw [style=new edge style 0] (95.center) to (93.center);
		\draw [style=new edge style 0] (86) to (97.center);
		\draw [style=new edge style 0] (99.center) to (88);
		\draw [style=new edge style 0] (89) to (101.center);
	\end{pgfonlayer}
\end{tikzpicture}}
\end{center}
\caption{An example of clustering-based network reduction. The original network (a) is clustered into (b) and the reduced network with the clusters is (c). The input and output structure after reduction is the same as the original network.}
\label{fig:clustering example}
\end{figure*}

As shown in \cite{Cheng2020ModelRM}, we can obtain a reduced model by aggregating clusters $\mathcal{C}_1, \mathcal{C}_2,\ldots ,\mathcal{C}_r$, that are illustrated in Fig. \ref{fig:clustering example}, into the nodes of $\mathcal{G}_r$ as follows.
\begin{align}
    \label{eq_clustering}
  \begin{cases}
    \Pi^\mathsf{T}\Pi\dot{x}_r (t) = \Pi^\mathsf{T} A \Pi x_r (t) + \Pi^\mathsf{T}B u(t), \\
    y_r (t) = C\Pi x_r (t), 
  \end{cases}
\end{align}
where $\Pi \in \mathbb{R}^{n\times r}$ is the characteristic matrix
and
\begin{align*}
    (\Pi)_{ij}=\left\{
    \begin{array}{ll}
        1, & \text{if } i \in \mathcal{C}_j,\\
        0, & \text{otherwise}.
    \end{array}
    \right. 
\end{align*}
In this case, reduced system \eqref{eq_clustering} has the same interconncetion structure with the original network system \eqref{original}.
That is,
if there is a directed path from $i\in \mathcal{V}_r$ to $j\in \mathcal{V}_r$, there is a directed path from 
a node of $\pi^{-1}(i)$ to a node of $\pi^{-1}(j)$ in $\mathcal{G}$.
Here, $\pi:\mathcal{V}\rightarrow \mathcal{V}_r$ is the associated map to the characteristic matrix $\Pi$ of the clustering $\{C_1,\ldots, C_r\}$.
We remark that $\Pi^\mathsf{T}\Pi$ is a diagonal matrix whose diagonal elements are the number of nodes in each cluster.

\subsection{Riemannian Manifold of Stable Matrices}

We will explain that the space of stable matrices $\mathbb{S}_r$ can be regarded as a Riemannian manifold.

As shown in Proposition 1 of \cite{prajna2002lmi},
for any stable matrix $A_r\in \mathbb{R}^{r\times r}$, there exists a $(J_r, R_r, Q_r)\in \mathbb{S}_r$ satisfying $A_r = (J_r-R_r)Q_r$. Conversely, for any $(J_r, R_r, Q_r)\in \mathbb{S}_r$, $(J_r-R_r)Q_r$ is stable. Here, $\mathbb{S}_r := \skew{r}\times \sym{r}\times \sym{r}$, $\skew{r}$ is the set of all skew symmetric matrices, and $\sym{r}$ is the set of all symmetric positive definite matrices.

The Euclidean space $\mathbb{R}^{k\times l}$ is a Riemannian manifold with the metric 
  $\left\langle \xi, \eta \right\rangle_P^{(\mathbb{R}^{k\times l})}\coloneqq \tr(\xi^{\mathsf{T}} \eta)$.
For $\sk{r}$, which is an Euclidean embedded submanifold endowed with the restricted metric
$\left\langle \xi, \eta \right\rangle_P^{(\mathbb{R}^{r\times r})}|_{\sk{r}}$ and also linear subspace of $\mathbb{R}^{r\times r}$,
the Riemannian gradient of $f:\sk{r}\to \mathbb{R}$ is calculated from the Euclidean gradient $\grad{}\bar{f}(x)$
easily using the orthogonal projection onto the tangent space; 
  $\rgrad f(x) = \mathrm{skew} (\grad{\bar{f}} (x)) \coloneqq \frac{1}{2}(\grad{\bar{f}} (x) - \grad{\bar{f}} (x)^{\mathsf{T}})$.
For more detail discussion, see \cite[Chapter~3]{boumal2020introduction} and \cite{absil2009optimization}.
The Riemannian metric of $\sym{r}$ is defined as 
  $\left\langle \xi, \eta \right\rangle_P^{(\sym{r})}\coloneqq \tr(P^{-1}\xi P^{-1}\eta)$,
for any $P\in\sym{r}$ and $\xi, \eta \in T_P\sym{r}$ \cite[Chapter~X\hspace{-.1em}I\hspace{-.1em}I]{lang2012fundamentals}.
The Riemannian manifold $\sym{r}$ with this metric has a closed form of the exponential map
  $\mathrm{Exp}_P({\xi}) =  P\exp(P^{-1}\xi)$, where $\exp(P)$ is a matrix exponential,
and the Riemannian gradient on $\sym{r}$ is calculated as
    $\rgrad f(P) =  P \mathrm{sym}(\grad{\bar{f}(P)}) P$
  by letting $\mathrm{sym}(S)\coloneqq \frac{1}{2}\qty(S + S^{\mathsf{T}})$ \cite{sato2017structure}.
  
  The set $\mathbb{S}_r$ is a product of the Riemannian manifolds $\skew{r}$ and $\sym{r}$. Thus, $\mathbb{S}_r$ is a Riemannian manifold with the canonically induced metric. 

\section{Problem formulation}

\subsection{Initial Reduced Network}

To define the constraints for the interconnection structure, we calculate the reduced model (\ref{eq:reduced system all}) on $\mathcal{G}_r=(\mathcal{V}_r, \mathcal{E}_r)$ by using the clustering-based algorithm as follows;
\begin{align}
\label{shokiten}
\begin{cases}
   A_r^{(0)} \coloneqq (\Pi^\mathsf{T}\Pi)^{-1}\Pi^\mathsf{T} A \Pi - \alpha I_r, \\
    B_r^{(0)} \coloneqq (\Pi^\mathsf{T}\Pi)^{-1}\Pi^\mathsf{T}B,\\
    C_r^{(0)} \coloneqq C\Pi
\end{cases}
 \end{align}
where $\alpha \geq 0$ is a sufficiently large constant such that $A_r^{(0)}$ is stable.
Note that
$(\Pi^\mathsf{T}\Pi)^{-1}\Pi^\mathsf{T} A \Pi$ is not stable in general, even if $A$ is stable.
For example, consider 
\begin{align*}
    A = \begin{bmatrix}
    -1 & 0 & 0 & 0 & 0\\
    6 & -2 & 0 & 0 & 0 \\
    0 & 1 & -3 & 0 & 0 \\
    0 & 0 & 1 & -4 & 0 \\
    0 & 0 & 0 & 10 & -5
    \end{bmatrix},\,\,
    \Pi = \begin{bmatrix}
    1 & 0 \\
    1 & 0 \\
    1 & 0 \\
    0 & 1 \\
    0 & 1
    \end{bmatrix}, 
\end{align*}
where $A$ is a stable and Metzler matrix.
Then, $
    (\Pi^\mathsf{T}\Pi)^{-1}\Pi^\mathsf{T} A \Pi = \begin{bmatrix}
    \frac{1}{3} & 0 \\
    \frac{1}{2} & \frac{1}{2}
    \end{bmatrix}$,
which is Metzler, but is not stable.
It is also notable that the non-diagonal nonzero entries of $A_r$ correspond to the edges of $\mathcal{G}_r$. Also, the nonzero entries of $B_r$ and $C_r$ imply the input and output position, respectively.

To preserve the interconnection structure, every entry of the feasible solution $(A_r, B_r, C_r)$ should be zero if the corresponding entry of the initial matrix $(A_r^{(0)}, B_r^{(0)}, C_r^{(0)})$ is zero and should be nonnegative if the corresponding entry is nonnegative.
For convenience, $\idz(X)$ denotes the indices sets of the zero-entries of the matrix $X$ and $\idnz(X)$ denotes the indices sets of the nonnegative-entries of the matrix $X$ for $X = A_r^{(0)}, B_r^{(0)}, C_r^{(0)}$. It is notable that $\idz(A_r^{(0)})$ and $\idnz(A_r^{(0)})$ do not contain the diagonal indices of $A_r^{(0)}$.

\subsection{$H^2$ Optimization Problem with constraints} \label{Sec2-D}

We
consider
an $H^2$ optimal model reduction problem using
the transfer functions $G$ of (\ref{original}) and $G_r$ of (\ref{eq:reduced system all}) 
to reconstruct a novel reduced model $(A_r, B_r, C_r)$ of preserving the positivity and interconnection structure better than a given reduced model $(A_r^{(0)}, B_r^{(0)}, C_r^{(0)})$ in the sense of the $H^2$ norm.
This is because
$\sup_{t\geq 0}\|y(t)-y_r(t)\|_2 \leq \|G-G_r\|_{H^2}$
assuming that $\|u\|_{L^2}\leq 1$, as explained in \cite{gugercin2008h_2, sato2019riemannian}.
This inequality indicates that
the maximum output error norm can be expected to become almost zero when
 $\|G-G_r\|^2_{H^2} = 2F(A_r, B_r, C_r) + \|G\|_{H^2}^2$ is sufficiently small,
 where
\begin{align}
    F(A_r, B_r, C_r) := & \frac{1}{2}\mathrm{tr}\qty(C_rSC_r^\mathsf{T} - 2C_rX^\mathsf{T} C^\mathsf{T}) \label{objective_function}\\
    = & \frac{1}{2}\mathrm{tr}\qty(B_r^\mathsf{T} SB_r + 2B^\mathsf{T} Y B_r ). \nonumber
\end{align} 
Here,
$(X, Y, S, T)$ are the solutions to the following Sylvester equations
\begin{align}
    \label{eq:4syls}
    AX + XA_r^\mathsf{T} + BB_r^\mathsf{T} &= 0,\\
    \label{eq:4sylsbig}
    A^\mathsf{T} Y + YA_r - C^\mathsf{T} C_r &=0,\\
    \label{eq:4sylsS}
    A_rS + SA_r^\mathsf{T} + B_rB_r^\mathsf{T} &=0,\\
    \label{eq:4syls2}
    A_r^\mathsf{T} T + TA_r + C_r^\mathsf{T} C_r &= 0.
  \end{align}

Therefore, 
we formulate the following optimization problem.
\begin{empheq}{align}
  \begin{aligned}
  &\underset{\substack{((J_r, R_r, Q_r), B_r, C_r)\\ \in {M_r}}}{\text{minimize}}&&F(A_r, B_r, C_r) \label{opt:formulated}\\
  &\text{subject to} &&\hspace{-10mm}  A_r \coloneqq \qty(J_r - R_r)Q_r,\\
  & &&\hspace{-10mm} g_{ij}^{({A_r})}(A_r)\leq 0, {}^\forall(i, j) \in \idnz(A_r^{(0)}),\\
  & &&\hspace{-10mm} g_{ij}^{({B_r})}(B_r)\leq 0, {}^\forall (i, j) \in \idnz(B_r^{(0)}),\\
  & &&\hspace{-10mm} g_{ij}^{({C_r})}(C_r)\leq 0, {}^\forall (i, j) \in \idnz(C_r^{(0)}),\\
  & &&\hspace{-10mm} g_{ij}^{({A_r})}(A_r)= 0, {}^\forall(i, j) \in \idz(A_r^{(0)}),\\
  & &&\hspace{-10mm} g_{ij}^{({B_r})}(B_r)= 0, {}^\forall (i, j) \in \idz(B_r^{(0)}),\\
  & &&\hspace{-10mm} g_{ij}^{({C_r})}(C_r)= 0, {}^\forall (i, j) \in \idz(C_r^{(0)}),\\
  \end{aligned}
\end{empheq}
where $M_r := \mathbb{S}_r \times \mathbb{R}^{r\times m} \times \mathbb{R}^{p\times r}$ and $g_{ij}^{({X})}(X):=- (X)_{ij}$ for $X=A_r, B_r, C_r$ and $(i, j)$ is an index of $X$.

The problem \eqref{opt:formulated} is an $H^2$ optimal model reduction problem with nonnegativity and interconnection structure-preserving constraints.
The reduced system (\ref{eq:reduced system all}) corresponding to a feasible solution to \eqref{opt:formulated} is always an asymptotically stable, positive, and has the same interconnection structure with the original system \eqref{original}.

Moreover, the optimization problem \eqref{opt:formulated} can be regarded as a Riemannian optimization problem with the constraints by introducing the Riemannian metric
\begin{multline}
    \langle\eta, \xi\rangle_x^{(M_r)} := 
\left\langle \xi, \eta \right\rangle_{J_r}^{(\mathbb{R}^{r\times r})}|_{\sk{r}}+    
    \left\langle \xi, \eta \right\rangle_{R_r}^{(\sym{r})} \\
    + 
    \left\langle \xi, \eta \right\rangle_{Q_r}^{(\sym{r})}
    +\left\langle \xi, \eta \right\rangle_{B_r}^{(\mathbb{R}^{r\times m})}
    +\left\langle \xi, \eta \right\rangle_{C_r}^{(\mathbb{R}^{p\times r})} \label{Mr_metric}
\end{multline}
into the set $M_r$,
where $x = (J_r, R_r, Q_r, B_r, C_r)$, and 
$\xi=(\xi_1, \xi_2, \xi_3, \xi_4, \xi_5)$, $\eta=(\eta_1, \eta_2, \eta_3, \eta_4, \eta_5)\in T_xM_r$. 
Therefore, we can develop an algorithm for solving \eqref{opt:formulated} based on Riemannian optimization \cite{absil2009optimization}.

\begin{remark}
According to Theorem 4.C.2 in \cite{takayama1985mathematical}, if the matrix $A$ has a dominant diagonal that is negative, $A$ is stable.
Using this fact, we can formulate another optimization problem by adding the inequality constraints to enforce the strict diagonally dominance.
However, this is just a sufficient condition for $A$ to be stable unlike our formulation in \eqref{opt:formulated}.
That is, our formulation is
useful to decrease the objective function value compared with the addition of the inequality constraints,
because the search space of $A_r$ is wider.
\end{remark}

\section{Proposed Method} \label{Sec3}



Because the optimiztion problem \eqref{opt:formulated} is a Riemannian optimization problem with constraints,
we develop an algorithm for solving \eqref{opt:formulated} based on RALM proposed in \cite{liu2019simple}.

The Lagrangian function of (\ref{opt:formulated}) for RALM is as follows:
\footnotesize
\begin{multline}
    \mathcal{L}_\rho ((J_r, R_r, Q_r, B_r,C_r), \lambda, \gamma)
  = F(A_r, B_r, C_r)\\
  + \frac{\rho}{2}\sum_{V\in \qty{A_r, B_r, C_r}}\Large\left( \sum _{(i, j)\in \idnz(V^{(0)})} \max \qty{0, \frac{\lambda_{ij}^{(V)}}{\rho} + g_{ij}^{(V)}(V)}^2\Large\right. \\
  \Large\left. +\sum _{(i, j)\in \idz(V^{(0)})} \qty({g_{ij}^{(V)}(V)} + \frac{\gamma^{V} _{ij}}{\rho})^2\Large\right), \nonumber
\end{multline}
\normalsize
where $\rho > 0$ is a penalty parameter and $\gamma_{ij}^{A_r} \in \mathbb{R}^{r\times r}$, $\gamma_{ij}^{B_r} \in \mathbb{R}^{r\times m}$, $\gamma_{ij}^{C_r} \in \mathbb{R}^{p\times r}$, $\lambda_{ij}^{A_r} \in \mathbb{R}^{r\times r}_{\geq \bm{0}}$, $\lambda_{ij}^{B_r} \in \mathbb{R}^{r\times m}_{\geq \bm{0}}$, and $\lambda_{ij}^{C_r} \in \mathbb{R}^{p\times r}_{\geq \bm{0}}$ are the hyper parameters of RALM.

The algorithm is shown in Algorithm \ref{alg:RALM}. Here, in Algorithm \ref{alg:RALM}, $\mathrm{dist}$ is the distance function on the Riemannian manifold $M_r$ equipped with the Riemannian metric \eqref{Mr_metric}.
That is,
\begin{multline}
    \dist{M_r} (x, y)^2 = \norm{J_r - J_r'}_{\mathrm{F}}^2 
    + \norm{\log R_r^{-1/2}R_r'R_r^{-1/2}}_{\mathrm{F}}^2 \\
     + \norm{\log Q_r^{-1/2}Q_r'Q_r^{-1/2}}_{\mathrm{F}}^2 
    + \norm{B_r - B_r'}_{\mathrm{F}}^2
    + \norm{C_r - C_r'}_{\mathrm{F}}^2\nonumber
\end{multline}
where $x = (J_r, R_r, Q_r, B_r, C_r)$, $y = (J_r', R_r', Q_r', B_r', C_r')\in M_r$ and 
$\xi=(\xi_1, \xi_2, \xi_3, \xi_4, \xi_5)$, $\eta=(\eta_1, \eta_2, \eta_3, \eta_4, \eta_5)\in T_xM_r$.
For the details of step 9 in Algorithm \ref{alg:RALM}, see \cite{absil2009optimization}.

To solve the subproblem of step 3 in Algorithm \ref{alg:RALM} by using a Riemannian line search method~\cite{absil2009optimization}, we calculate the Euclidean gradients of $\lag$.
As shown in \cite{van2008h2}, the Euclidean gradients $\bar{F}$ with respect to $A_r, B_r$ and $C_r$
are calculated as
\begin{align}
  \grad{} _{A_r} \bar{F} = TS + Y^\mathsf{T} X,\, \,
  &\grad{} _{B_r} \bar{F} = TB_r + Y^\mathsf{T} B,\nonumber\\
  \grad{} _{C_r} \bar{F} = &C_rS - CX,\nonumber
\end{align}
respectively,
where $(X,\ Y,\ S,\ T)$ are the solutions to (\ref{eq:4syls})-(\ref{eq:4syls2}). 
Besides, it is easily seen that
\begin{align}
 \sum_{(i, j)\in \idnz ({A_r^{(0)}})} \grad{} _{A_r} \bar{g}_{ij}^{(A_r)} &= \qty(I_r - \mathds{1}(r, r))\odot \chi_{\idnz(A_r^{(0)})},\nonumber\\
  \sum_{(i, j)\in \idnz ({B_r^{(0)}})} \grad{} _{B_r} \bar{g}_{ij}^{(B_r)} &= -\chi_{\idnz(B_r^{(0)})},\nonumber\\
  \sum_{(i, j)\in \idnz ({C_r^{(0)}})} \grad{} _{C_r} \bar{g}_{ij}^{(C_r)} &= - \chi_{\idnz(C_r^{(0)})},\nonumber\\
  \sum_{(i, j)\in \idz ({A_r^{(0)}})} \grad{} _{A_r} \bar{g}_{ij}^{(A_r)} &= \qty(I_r - \mathds{1}(r, r))\odot \chi_{\idz(A_r^{(0)})},\nonumber\\
  \sum_{(i, j)\in \idz ({B_r^{(0)}})} \grad{} _{B_r} \bar{g}_{ij}^{(B_r)} &= -\chi_{\idz(B_r^{(0)})},\nonumber\\
  \sum_{(i, j)\in \idz ({C_r^{(0)}})} \grad{} _{C_r} \bar{g}_{ij}^{(C_r)} &= - \chi_{\idz(C_r^{(0)})},\nonumber
\end{align}
where $\odot $ denotes element-wise product and $\mathds{1}(k, l)$ is a $k\times l$ matrix whose entries are all $1$.
Here, $\chi_{\idnz{(V)}}$ and $\chi_{\idz{(V)}}$ are matrices of the same shape as $(V)$, being defined as 
\begin{align}
    \qty(\chi_{\idnz{(V)}})_{i, j} &= \begin{cases}
        1, & \text{if } (i,j)\in {\idnz{(V)}},\\
        0, & \text{otherwise},
    \end{cases}\nonumber\\
    \qty(\chi_{\idz{(V)}})_{i, j} &= \begin{cases}
        1, & \text{if } (i,j)\in {\idz{(V)}},\\
        0, & \text{otherwise},
    \end{cases}\nonumber
\end{align}
for $V\in \qty{A_r^{(0)}, B_r^{(0)}, C_r^{(0)}}$.
Then, the Euclidean gradient of the Lagrangian is written as
\begin{align}
  \grad{} _{A_r} &\overline{\mathcal{L}}_{\rho} ((J_r, R_r, Q_r,B_r,C_r), \lambda, \gamma)\nonumber \\
  = &\ \grad{} _{A_r} \bar{F}(A_r, B_r, C_r)\nonumber \\
  &+ \Bigl( U_{\mathrm{i}}^{(A_r)}\odot \chi_{\idnz(A_r^{(0)})}+ U_{\mathrm{e}}^{(A_r)}\odot \chi_{\idz(A_r^{(0)})}\Bigr)\nonumber \\ 
  &\phantom{+ \Bigl(}\odot (I_r - \mathds{1}(r, r)),\nonumber\\
  \grad{} _{B_r} &\overline{\mathcal{L}}_{\rho} ((J_r, R_r, Q_r,B_r,C_r), \lambda, \gamma)\nonumber \\
  = &\grad{} _{B_r} \bar{F}(A_r, B_r, C_r)\nonumber \\
  & -  \Bigl( U_{\mathrm{i}}^{(B_r)}\odot \chi_{\idnz(B_r^{(0)})} + U_{\mathrm{e}}^{(B_r)}\odot \chi_{\idz(B_r^{(0)})}\Bigr),\nonumber\\
  \grad{} _{C_r} &\overline{\mathcal{L}}_{\rho} ((J_r, R_r, Q_r,B_r,C_r), \lambda, \gamma)\nonumber \\
  =& \grad{} _{B_r} \bar{F}(A_r, B_r, C_r)\nonumber \\
  & - \Bigl( U_{\mathrm{i}}^{(C_r)}\odot \chi_{\idnz(C_r^{(0)})} + U_{\mathrm{e}}^{(C_r)}\odot \chi_{\idz(C_r^{(0)})}\Bigr),\nonumber
\end{align}
and
\begin{align}
    U_{\mathrm{i}}^{(V)}&\coloneqq \qty(\max \qty{0, \frac{\lambda_{ij}^{(V)}}{\rho} + g_{ij}^{(V)}(V)})_{(i, j)\in \mathcal{I}_V}\nonumber\\
    U_{\mathrm{e}}^{(V)}&\coloneqq \qty(\frac{\lambda_{ij}^{(V)}}{\rho} + g_{ij}^{(V)}(V))_{(i, j)\in \mathcal{I}_V}\nonumber
\end{align}
 are matrices for each $V \in \qty{A_r, B_r, C_r}$.
 Using the chain rule,
 we obtain the Euclidean gradients
  $\grad{} _{J_r} \blag = (\grad{} _{A_r} \blag ) Q_r$, $\grad{} _{R_r} \blag
  = - (\grad{} _{A_r} \blag ) Q_r$, and
  $\grad{} _{Q_r} \blag
  = -(J_r+R_r) (\grad{} _{A_r} \blag )$.

Finally, we calculate the Riemannian gradients from the Euclidean gradients in the same manner described in Section~II-D.
The Riemannian gradients are used to solve the subproblem of step 3 in Algorithm \ref{alg:RALM}.

\floatstyle{spaceruled}
\restylefloat{algorithm}
\begin{algorithm}
\caption{Riemannian augmented Lagrangian method (RALM).}
\label{alg:RALM}
     \begin{algorithmic}[1]
     \renewcommand{\algorithmicensure}{\textbf{Input:}}
      \Require The function $F$ in \eqref{objective_function} and the constraints $\qty{g_{ij}^{(V)}}_{(i, j) \in \mathcal{I}_V, V\in \qty{A_r, B_r, C_r}}$ in \eqref{opt:formulated} on Riemannian manifold $M_r$.
      \Ensure Initial point $x_0 = (J^{(0)}, R^{(0)}, Q^{(0)}, B^{(0)}, C^{(0)})\in M_r$;
    initial hyper parametars $\sigma = (\sigma_{A_r}, \sigma_{B_r}, \sigma_{C_r}) \in \mathbb{R}^{r\times r}\times \mathbb{R}^{r\times m}\times \mathbb{R}^{p\times r}$, $\rho\in \mathbb{R}_{\geq 0}$, $\lambda = (\lambda_{A_r}, \lambda_{B_r}, \lambda_{C_r}) \in \mathbb{R}_{\geq 0}^{r\times r}\times \mathbb{R}_{\geq 0}^{r\times m}\times \mathbb{R}_{\geq 0}^{p\times r}$, and $\gamma = (\gamma_{A_r}, \gamma_{B_r}, \gamma_{C_r}) \in \mathbb{R}^{r\times r}\times \mathbb{R}^{r\times m}\times \mathbb{R}^{p\times r}$;
    initial accuracy tolerance $\varepsilon > 0$; minimum tolerance $\varepsilon _\mathrm{min} > 0$ s.t. $\varepsilon_\mathrm{min} < \varepsilon$;
    constants $\theta_\rho>1, \theta_\varepsilon\in\qty(0, 1); \theta_\sigma\in \qty(0, 1)$;
    boundary vectors $\lambda^{\mathrm{min}}=(\lambda_{A_r}^{\mathrm{min}}, \lambda_{B_r}^{\mathrm{min}}, \lambda_{C_r}^{\mathrm{min}}), \lambda^{\mathrm{max}}=(\lambda_{A_r}^{\mathrm{max}}, \lambda_{B_r}^{\mathrm{max}}, \lambda_{C_r}^{\mathrm{max}})\in \mathbb{R}_{\geq 0}^{r\times r}\times \mathbb{R}_{\geq 0}^{r\times m}\times \mathbb{R}_{\geq 0}^{p\times r}$ s.t.
    $\lambda_{V}^{\mathrm{min}}\leq \lambda_{V} \leq \lambda_{V}^{\mathrm{max}}$, and $\gamma^{\mathrm{min}}=(\gamma_{A_r}^{\mathrm{min}}, \gamma_{B_r}^{\mathrm{min}}, \gamma_{C_r}^{\mathrm{min}}), \gamma^{\mathrm{max}}=(\gamma_{A_r}^{\mathrm{max}}, \gamma_{B_r}^{\mathrm{max}}, \gamma_{C_r}^{\mathrm{max}})\in \mathbb{R}^{r\times r}\times \mathbb{R}^{r\times m}\times \mathbb{R}^{p\times r}$ s.t.
    $\gamma_{V}^{\mathrm{min}}\leq \gamma_{V} \leq \gamma_{V}^{\mathrm{max}}$ for $V\in \qty{A_r, B_r, C_r}$;
    minimum step size $d_{\mathrm{min}}>0$.
      \For{$k=0, 1,\cdots $}
        \State Solve the subproblem below using Riemannian line search methods within a tolerance $\varepsilon$:
        \State\[
        \underset{x\in M_r}{\text{min}}\quad \lagr{{}}(x, \lambda, \gamma)
        \]
      \State and set the solution as $x_{k+1} = (J^{(k+1)}, R^{(k+1)}, Q^{(k+1)}, B^{(k+1)}, C^{(k+1)})$\;
      \State Set $A^{(k+1)} = \qty(J^{(k+1)} - R^{(k+1)})Q^{(k+1)}$\;
      \If{$\mathrm{dist}^{(M_r)}(x_k, x_{k+1})<d_{\mathrm{min}}$ and $\varepsilon \leq \varepsilon_{\mathrm{min}}$}
        \State  \Return{$x_{k+1}$}\;
      \EndIf
      \State Update the hyper parameters $\lambda$, $\gamma$, $\sigma$, $\rho$, and accuracy tolerance $\varepsilon$.
    \EndFor
    \end{algorithmic}
\end{algorithm}

\begin{remark}
For the time complexity of the proposed method, the bottleneck of the algorithm is to calculate the objective function $F$ and its gradients because $F$ internally requires the solutions for large-scale Sylvester equations \eqref{eq:4syls} and \eqref{eq:4sylsbig}.
However, in many applications, these Sylvester equations have the sparse-dense structure. That is, 
\begin{enumerate}
    \item the original matrix $A$ is large-scale, but is sparse.
    \item the reduced matrix $A_r$ is small-scale, but is dense.
\end{enumerate}
In this situation, we can use an efficient algorithm whose computational complexity is greatly smaller than $O(n^3)$ for solving  \eqref{eq:4syls} and \eqref{eq:4sylsbig} such as the method proposed in Section 3 in \cite{benner2011sparse}.
\end{remark}

\section{Experiment}
\subsection{Experimental Conditions}
\begin{table*}[t]
\centering
\caption{RALM parameters used in experiment}
  \label{table:param}
\begin{tabular}{l||llllllllllllll}
parameter & $\sigma^{(0)}$ & $\rho_0$ & $\lambda^{(0)}$ & $\lambda^{\mathrm{min}}$ & $\lambda^{\mathrm{max}}$ & $\gamma^{(0)}$ & $\gamma^{\mathrm{min}}$ & $\gamma^{\mathrm{max}}$ & $\varepsilon_0$ & $\varepsilon_{\mathrm{min}}$ & $\theta_\rho$ & $\theta_\varepsilon$ & $\theta_\sigma$ & $d_\mathrm{min}$           \\ \hline
value     & 0              & 10       & 1.0             & 0                        & 10                       & 0              & -1.5                    & 1.5                     & 1.0             & $10^{-16}$   & 1.01        & 0.95                 & 0.9             & $10^{-8}$
\end{tabular}
\end{table*}
We conducted a numerical experiment to verify the effectiveness of the proposed method in the sense of $H^2$-error and $H^\infty$-error.
We used the network shown in Fig. \ref{fig:clustering example} with the random positive weight sampled independently from the uniform distribution on $[0, 1]$. Then, the coefficient matrices are
\begin{align}
    A &= - \mathcal{L} - 0.1I_n,\nonumber\\
    B_{ij} &= \begin{cases}
       1 & \text{if }(i, j)=(1, 1) \text{ or } (14, 2)\\
       0 & \text{otherwise.}
    \end{cases},\\
    C_{ij} &= \begin{cases}
       1 & \text{if }(i, j)=(1, 1) \text{ or } (2, 7) \text{ or } (3, 18)\\
       0 & \text{otherwise.}
    \end{cases},\nonumber
\end{align}
where $\mathcal{L}$ is a loopy Laplacian of the random weighted network and the second term of $A$ is for numerical stability.
We used Riemannian conjugate gradient descent method \cite{absil2009optimization} with the Riemannian gradients obtained in Section III  as the subsolver in Algorithm~\ref{alg:RALM}.
For each iteration, we used the after-100-iteration output of the Riemannian conjugate gradient descent method as its solution.
The parameters for Algorithm~\ref{alg:RALM} are shown in Table~\ref{table:param}.

We obtained the initial iterative point $(J_r^{(0)}, R_r^{(0)}, Q_r^{(0)}, B_r^{(0)}, C_r^{(0)})$ by the following way.
\begin{enumerate}
    \item Calculate $A_r^{(0)}$, $B_r^{(0)}$, and $C_r^{(0)}$ using (\ref{shokiten}).
    \item Solve the Lyapunov equation ${A_r^{(0)}}^\mathsf{T}Q_r^{(0)} + Q_r^{(0)}A_r^{(0)} = -I_r$ for $Q_r^{(0)}$.
    \item Calculate 
        $J_r^{(0)} =\frac{1}{2}(A_r^{(0)}(Q_r^{(0)})^{-1} - (Q_r^{(0)})^{-1}{A_r^{(0)}}^\mathsf{T})$ and
        $R_r^{(0)} = -\frac{1}{2}({A_r^{(0)}}(Q_r^{(0)})^{-1} + (Q_r^{(0)})^{-1}{A_r^{(0)}}^\mathsf{T})$.
\end{enumerate}
Here, $(J_r^{(0)}, R_r^{(0)}, Q_r^{(0)})\in \mathbb{S}_r$, as shown in \cite{sato2019riemannian}.

We define
the relative $H^2$ and $H^\infty$ errors as
\begin{align}
    \mathrm{Err}_{2}(G_r) :=\frac{\|G - G_r\|_{H^2}}{\|G\|_{H^2}},\,
    \mathrm{Err}_{\infty}(G_r) :=\frac{\|G - G_r\|_{H^\infty}}{\|G\|_{H^{\infty}}}.\nonumber
\end{align}



\subsection{Result}
\floatstyle{plain}
\restylefloat{figure}
\begin{figure}[t]
\vspace{2mm}
\begin{center}
\includegraphics[width=75mm]{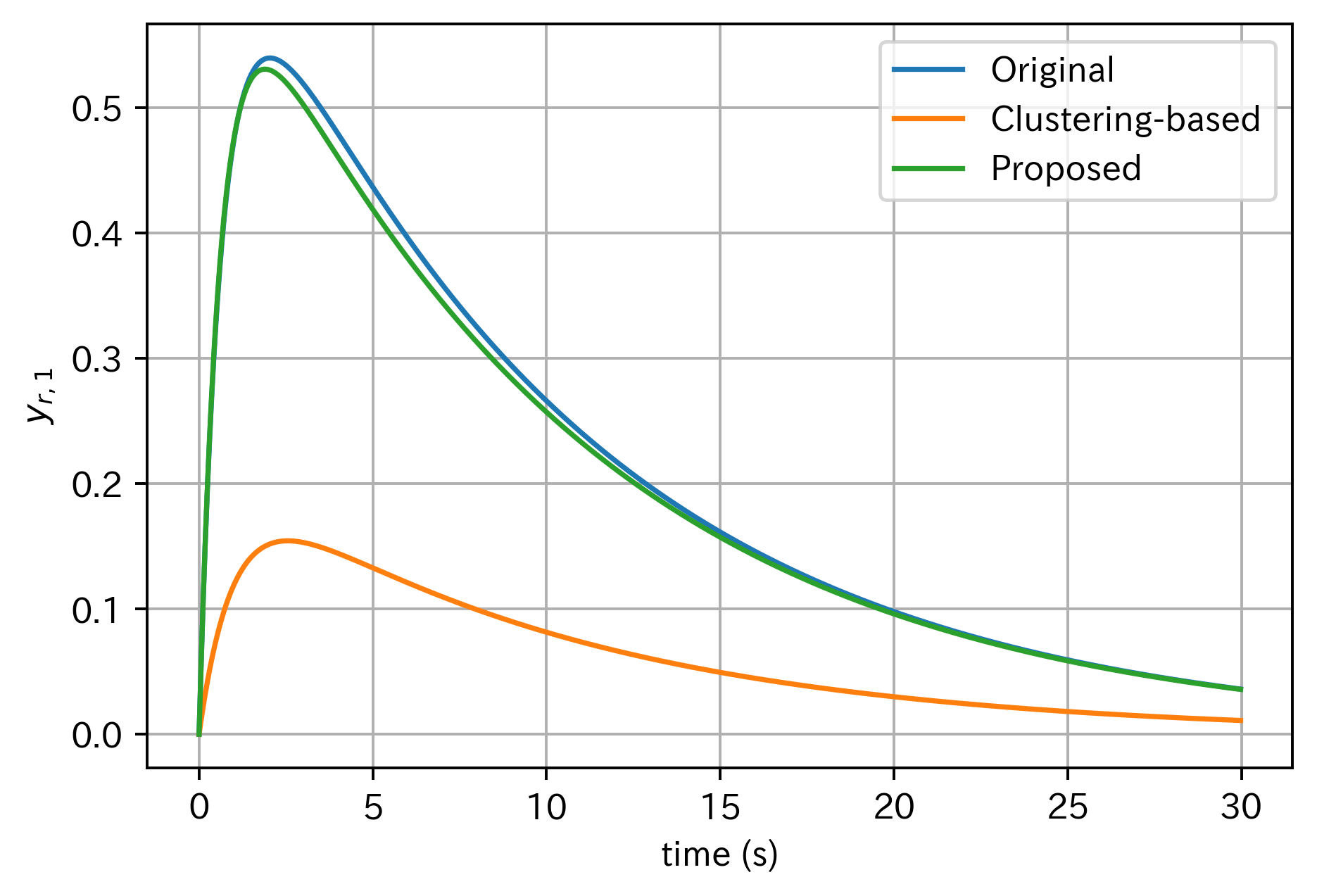}
\end{center}
\caption{The output $y_{r, 1}(t)$ for each original system (blue), reduced model by clustering-based method (orange), and reduced model by proposed method (green)}
\label{fig:output example}
\end{figure}
The $H^2$-error of clustering-based method was $\mathrm{Err}_2(G_r)=75.7\%$ and its $H^\infty$-error was $\mathrm{Err}_\infty(G_r)=70.24\%$. On the other hand, the $H^2$-error of propsosed method was $\mathrm{Err}_2(G_r)=3.14\%$ and its $H^\infty$-error was $\mathrm{Err}_\infty(G_r)=4.67\%$.
It is easily seen that the proposed method improves the clustering-based method in the sense of not only $H^2$-error but also $H^\infty$-error. Moreover, in Fig.~\ref{fig:output example}, we illustrate the example outputs corresponding to the input
\begin{align}
    u(t)=\begin{pmatrix}
    \exp(-0.1t)|\cos (100\pi t)|\\  \exp(-0.1t)|\sin (200\pi t)|\end{pmatrix}.\nonumber
\end{align}

\section{Conclusion}

We formulated the model reduction problem of an asymptotically stable and positive network system as a constrained Riemannian optimization problem with $H^2$-error of the original and reduced network systems as the objective function.
Our method reduces the dimension of the network system while preserving the properties of stability, positivity, and interconnection structure by applying RALM and
deriving the Riemannian gradients of the Lagrangian.
We proposed to use the initial point of the clustering-based method.
We conducted a numerical experiment and compare it with the clustering-based method in the sense of $H^2$-error and $H^\infty$-error and verified that our method improved the reduction performance of the clustering-based method.

We note that our proposed algorithm can be easily extended to the case of a positive network system with multi-dimensional subsystems.
Moreover, instead of an initial model generated by the clustering method as explained in Section II-C, we can use other arbitrary reduction methods 
which preserve stability, positivity, and interconnection structure, as the initial model of our proposed algorithm.


%

\section*{Acknowledgment}
This work was supported by Japan
Society for the Promotion of Science KAKENHI under Grant 20K14760.

\bibliographystyle{IEEEtran}
\bibliography{IEEEabrv,under_grad_revised}

\begin{thebibliography}{10}
\providecommand{\url}[1]{#1}
\csname url@samestyle\endcsname
\providecommand{\newblock}{\relax}
\providecommand{\bibinfo}[2]{#2}
\providecommand{\BIBentrySTDinterwordspacing}{\spaceskip=0pt\relax}
\providecommand{\BIBentryALTinterwordstretchfactor}{4}
\providecommand{\BIBentryALTinterwordspacing}{\spaceskip=\fontdimen2\font plus
\BIBentryALTinterwordstretchfactor\fontdimen3\font minus
  \fontdimen4\font\relax}
\providecommand{\BIBforeignlanguage}[2]{{%
\expandafter\ifx\csname l@#1\endcsname\relax
\typeout{** WARNING: IEEEtran.bst: No hyphenation pattern has been}%
\typeout{** loaded for the language `#1'. Using the pattern for}%
\typeout{** the default language instead.}%
\else
\language=\csname l@#1\endcsname
\fi
#2}}
\providecommand{\BIBdecl}{\relax}
\BIBdecl

\bibitem{moore1981principal}
B.~Moore, ``Principal component analysis in linear systems: Controllability,
  observability, and model reduction,'' \emph{IEEE transactions on automatic
  control}, vol.~26, no.~1, pp. 17--32, 1981.

\bibitem{gugercin2008h_2}
S.~Gugercin, A.~C. Antoulas, and C.~Beattie, ``{$\mathcal{H}_2$} model
  reduction for large-scale linear dynamical systems,'' \emph{SIAM journal on
  matrix analysis and applications}, vol.~30, no.~2, pp. 609--638, 2008.

\bibitem{kellems2009low}
A.~R. Kellems, D.~Roos, N.~Xiao, and S.~J. Cox, ``Low-dimensional,
  morphologically accurate models of subthreshold membrane potential,''
  \emph{Journal of computational neuroscience}, vol.~27, no.~2, p. 161, 2009.

\bibitem{gugercin2012structure}
S.~Gugercin, R.~V. Polyuga, C.~Beattie, and A.~Van Der~Schaft,
  ``Structure-preserving tangential interpolation for model reduction of
  port-hamiltonian systems,'' \emph{Automatica}, vol.~48, no.~9, pp.
  1963--1974, 2012.

\bibitem{baur2011interpolatory}
U.~Baur, C.~Beattie, P.~Benner, and S.~Gugercin, ``Interpolatory projection
  methods for parameterized model reduction,'' \emph{SIAM Journal on Scientific
  Computing}, vol.~33, no.~5, pp. 2489--2518, 2011.

\bibitem{benner2015survey}
P.~Benner, S.~Gugercin, and K.~Willcox, ``A survey of projection-based model
  reduction methods for parametric dynamical systems,'' \emph{SIAM review},
  vol.~57, no.~4, pp. 483--531, 2015.

\bibitem{summers2014optimal}
T.~H. Summers and J.~Lygeros, ``Optimal sensor and actuator placement in
  complex dynamical networks,'' \emph{IFAC Proceedings Volumes}, vol.~47,
  no.~3, pp. 3784--3789, 2014.

\bibitem{gates2016control}
A.~J. Gates and L.~M. Rocha, ``Control of complex networks requires both
  structure and dynamics,'' \emph{Scientific reports}, vol.~6, no.~1, pp.
  1--11, 2016.

\bibitem{kim2018role}
J.~Z. Kim, J.~M. Soffer, A.~E. Kahn, J.~M. Vettel, F.~Pasqualetti, and D.~S.
  Bassett, ``Role of graph architecture in controlling dynamical networks with
  applications to neural systems,'' \emph{Nature physics}, vol.~14, no.~1, pp.
  91--98, 2018.

\bibitem{ishizaki2018graph}
T.~Ishizaki, A.~Chakrabortty, and J.-I. Imura, ``Graph-theoretic analysis of
  power systems,'' \emph{Proceedings of the IEEE}, vol. 106, no.~5, pp.
  931--952, 2018.

\bibitem{martin2018large}
N.~Martin, P.~Frasca, and C.~Canudas-de Wit, ``Large-scale network reduction
  towards scale-free structure,'' \emph{IEEE Transactions on Network Science
  and Engineering}, vol.~6, no.~4, pp. 711--723, 2018.

\bibitem{vandendorpe2008model}
A.~Vandendorpe and P.~Van~Dooren, ``Model reduction of interconnected
  systems,'' in \emph{Model order reduction: theory, research aspects and
  applications}.\hskip 1em plus 0.5em minus 0.4em\relax Springer, 2008, pp.
  305--321.

\bibitem{haddad2010nonnegative}
W.~M. Haddad, V.~Chellaboina, and Q.~Hui, \emph{Nonnegative and compartmental
  dynamical systems}.\hskip 1em plus 0.5em minus 0.4em\relax Princeton
  University Press, 2010.

\bibitem{jacquez1993qualitative}
J.~A. Jacquez and C.~P. Simon, ``Qualitative theory of compartmental systems,''
  \emph{Siam Review}, vol.~35, no.~1, pp. 43--79, 1993.

\bibitem{Cheng2020ModelRM}
X.~Cheng and J.~M.~A. Scherpen, ``Model reduction methods for complex network
  systems,'' \emph{Annual Review of Control, Robotics, and Autonomous Systems},
  vol.~4, pp. 425--453, 2020.

\bibitem{ISHIZAKI2015238}
T.~Ishizaki, K.~Kashima, A.~Girard, J.~ichi Imura, L.~Chen, and K.~Aihara,
  ``Clustered model reduction of positive directed networks,''
  \emph{Automatica}, vol.~59, pp. 238--247, 2015.

\bibitem{reis2009positivity}
T.~Reis and E.~Virnik, ``Positivity preserving model reduction,'' in
  \emph{Positive Systems}.\hskip 1em plus 0.5em minus 0.4em\relax Springer,
  2009, pp. 131--139.

\bibitem{sato2015riemannian}
H.~Sato and K.~Sato, ``{Riemannian trust-region methods for $H^2$ optimal model
  reduction},'' in \emph{54th IEEE Conference on Decision and Control (CDC)},
  2015, pp. 4648--4655.

\bibitem{sato2017structure}
K.~Sato and H.~Sato, ``Structure-preserving {$H^2$} optimal model reduction
  based on the riemannian trust-region method,'' \emph{IEEE Transactions on
  Automatic Control}, vol.~63, no.~2, pp. 505--512, 2018.

\bibitem{sato2018riemannian}
K.~Sato, ``{Riemannian optimal model reduction of linear port-Hamiltonian
  systems},'' \emph{Automatica}, vol.~93, pp. 428--434, 2018.

\bibitem{sato2019riemannian}
K.~Sato, ``Riemannian optimal model reduction of stable linear systems,''
  \emph{IEEE Access}, vol.~7, pp. 14\,689--14\,698, 2019.

\bibitem{absil2009optimization}
P.-A. Absil, R.~Mahony, and R.~Sepulchre, \emph{Optimization algorithms on
  matrix manifolds}.\hskip 1em plus 0.5em minus 0.4em\relax Princeton
  University Press, 2008.

\bibitem{cheng2019reduced}
X.~Cheng, L.~Yu, and J.~M. Scherpen, ``Reduced order modeling of linear
  consensus networks using weight assignments,'' in \emph{2019 18th European
  Control Conference (ECC)}, 2019, pp. 2005--2010.

\bibitem{liu2019simple}
C.~Liu and N.~Boumal, ``Simple algorithms for optimization on riemannian
  manifolds with constraints,'' \emph{Applied Mathematics \& Optimization}, pp.
  1--33, 2019.

\bibitem{prajna2002lmi}
S.~Prajna, A.~van~der Schaft, and G.~Meinsma, ``An {LMI} approach to
  stabilization of linear port-controlled {Hamiltonian} systems,''
  \emph{Systems \& control letters}, vol.~45, no.~5, pp. 371--385, 2002.

\bibitem{boumal2020introduction}
N.~Boumal, ``An introduction to optimization on smooth manifolds,''
  \emph{Available online, May}, 2020.

\bibitem{lang2012fundamentals}
S.~Lang, \emph{Fundamentals of differential geometry}.\hskip 1em plus 0.5em
  minus 0.4em\relax Springer Science \& Business Media, 2012, vol. 191.

\bibitem{takayama1985mathematical}
A.~Takayama, \emph{Mathematical economics}.\hskip 1em plus 0.5em minus
  0.4em\relax Cambridge university press, 1985.

\bibitem{van2008h2}
P.~Van~Dooren, K.~A. Gallivan, and P.-A. Absil, ``{$H_2$}-optimal model
  reduction of {MIMO} systems,'' \emph{Applied Mathematics Letters}, vol.~21,
  no.~12, pp. 1267--1273, 2008.

\bibitem{benner2011sparse}
P.~Benner, M.~K\"{o}hler, and J.~Saak, ``{Sparse-dense Sylvester equations in
  $H_2$-model order reduction},'' \emph{Max Planck Institute Magdeburg
  Preprints}, 2011.

\end{thebibliography}
\end{document}